\documentclass[a4paper,oneside]{article}
\RequirePackage[latin1]{inputenc}
\usepackage[french]{babel} 
\usepackage{amsmath} 
\usepackage{amssymb} 
\usepackage{graphicx} 
\usepackage[T1]{fontenc} 
\usepackage{times} 
\usepackage{eurosym} 
\usepackage{color} 
\usepackage{multicol} 
\usepackage{makeidx}
\usepackage{hyperref} 
\hypersetup{backref,colorlinks=true,bookmarks=true,pdftoolbar=true}
\newenvironment{preuve}[1][Proof :]{\textbf{#1} }{\hfill $\blacksquare$}
\numberwithin{equation}{section}
\newtheorem{theo}{Theorem}

\newtheorem{lem}{Lemma}
\newtheorem{cor}{Corollary}
\newtheorem{rem}{Remark}

\title{\bf{Some improvement on non-parametric estimation of income distribution and poverty index}}
\author{\bf{Youssou Ciss}\footnote{ciss.youssouf@gmail.com}\\
{LERSTAD, UFR SAT, Université Gaston Berger, B.P 234, Saint-Louis, Sénégal}\\
\\ \bf{El hadji Deme}\footnote{elhadjidemeufrsat@gmail.com}\\
{LERSTAD, UFR SAT, Université Gaston Berger B.P 234, Saint-Louis, Sénégal}\\
\\ \bf{Hamza Dhaker}\footnote{hamza.dhaker@umoncton.ca}\\
{D\'epartement de Math\'ematiques et Statistique, Universit\'e de Moncton, NB, Canada}}
\date{}
\begin{document}
\maketitle
\section*{Abstract}
\noindent\rule[2pt]{\textwidth}{0.5pt}
\noindent In this paper, we propose an estimator of Foster, Greer and Thorbecke class of measures $\displaystyle P(z,\alpha) = \int_0^{z}\Big(\frac{z-x}{z}\Big)^{\alpha}f(x)\, dx$, where $z>0$ is the poverty line, $f$ is the probabily density function of the income distribution and $\alpha$ is the so-called poverty aversion. The estimator is constructed with a bias reduced kernel estimator. Uniform almost sure consistency and uniform mean square consistenty  are established.
A simulation study indicates that our new estimator performs well.\\
\noindent\rule[2pt]{\textwidth}{0.5pt}

\noindent \textbf{Key words and expressions :} poverty line, bias reduction kernel, uniform almost sure consistency, uniforme  mean square consistency, rate of convergence.\\
\section{Introduction and definition of the estimator}
\noindent Let  \, $F$ \, be the cumulative distribution function of the income variable \, $X$ \, from a population with continuous density function \, $f$. The FGT (Foster, Greer, Thorbecke) \cite{Foster&Greer&Thorbecke1984} class of poverty index measures by the real \, $\alpha \geq 0$ \, is  defined by:
\begin{equation}
\begin{split}
\label{cc1}
&P(z,\alpha) =  \left\{\begin{array}{ccc}
\displaystyle \int_0^{z}\left(\frac{z-x}{z}\right)^{\alpha}f(x)\, dx \qquad &\text{if} \, z > 0, \\~~\\
0 \qquad &\text{otherwise}
\end{array}\right.
\end{split}
\end{equation}
where \, $z$ \, is the poverty line.\\

\noindent Let  \, $\displaystyle X_1, \cdot \cdot \cdot, X_n$ \, be a  random sample of size \, $n$ \, from income random variables (r.v.) with distribution function $F$. 
Seidl (1988) \cite{seidl1988} introduced
An empirical estimator of the FGT poverty index $P(z,\alpha)$ as following:
\[\widehat{P}_n(z,\alpha)=\frac{1}{n}\sum_{i=1}^n\Big(1-\frac{X_i}{z}\Big)_+^{\alpha} \qquad \text{where} \qquad x_+=\max(0, x). \quad \]
This estimator was fully useful in a large range of applications in economics (Widely used in practice in econometrics and actuarial). It is an unbiased consistent estimator and is asymptoticaly  normal. 
Lo {\em et al.} \cite{loal2009} used  empirical processes and extreme-values methodology to study this estimator. Seck \cite{ctseck2011}, Seck and Lo \cite{secka} used some non-weighted poverty measures, viewed as stochastic processes and indexed by real numbers or monotone functions, to follow up the poverty evolution between two periods.
Dia \cite{Galaye2008} and also Ciss {\em et al.} \cite{ciss1}  proposed  new kernels estimators, based on the Riemann sum, respectively, for $\alpha=0$ and $\alpha\geq 1$ and $\alpha \in ]0,1[$. 

\bigskip
\noindent Dia \cite{Galaye2008} and Ciss {\em et al.} \cite{ciss1} considered also the classical nonparametric estimator of the density  $f$  (Parzen-Rosenblatt): 
\[\hat{f}(x) = \frac{1}{n}\sum_{i=1}^n\frac{1}{h}K\Big(\frac{x-X_i}{h}\Big),\]
where  $h$ is a function of  $n$ which tends to zero as  $n$  tends to infinity and  $K$ verifies the following hypotheses:
\begin{equation} 
\label{cc2}(\mathbf{H_1}) \sup_{-\infty<x<+\infty}|K(x)|<+\infty,\quad
(\mathbf{H_2}) \int_{-\infty}^{+\infty}K(x)\, dx=1,\quad  
(\mathbf{H_3}) \lim_{x \rightarrow \pm\infty}|xK(x)|=0
\end{equation}
and proposed as  estimator of  FGT poverty index, the following one:
\begin{equation} 
\label{cc5}{P}_n(z,\alpha)=\frac{1}{n}\sum_{j=1}^n\sum_{i=0}^{[z/h]}\Big(\frac{z-ih}{z}\Big)^{\alpha}K\Big(\frac{ih-X_j}{h}\Big)\, .
\end{equation}

\bigskip
\noindent 
Recently, Zakaria et al. \cite{ciss5} considered the following adaptive kernel estimator of the density  $f$: 
\[\hat{f}^{\lambda}(x) = \frac{1}{n}\sum_{i=1}^n\frac{1}{h\lambda_i}K\Big(\frac{x-X_i}{h\lambda_i}\Big)\]
and proposed as  estimator of  FGT poverty index, the following estimator :
\begin{equation} 
\label{cc5}{P}^{\lambda}_n(z,\alpha)=\frac{1}{n}\sum_{j=1}^n\sum_{i=0}^{[z/h\lambda_j]}\Big(\frac{z-ih\lambda_j}{z}\Big)^{\alpha}K\Big(\frac{ih\lambda_j-X_j}{h\lambda_j}\Big)\,  \quad \text{for} \quad \alpha=0 \quad \text{or} \quad \alpha \geq 1,
\end{equation}
where $\lambda_j$,  is a parameter that varies according to the local concentration of the data. 

\bigskip
\noindent Both \cite{Parzen1962} and \cite{Silverman1986} pointed out that if $\int uK(u)du =0$ and $f$ is twice continuously differentiable in a neighborhood of $x$, then
\begin{equation}
\label{bias}
Bias(\hat{f}(x)) = E(\hat{f}(x)) - f(x) = \frac{1}{2}f''(x)h^2\int u^2K(u)du + O(h^3).
\end{equation}
In order to reduce the bias of the classical kernel estimator $\hat{f}$, we consider the following estimator introduced in \cite{Xie2014}
\begin{equation*}
\tilde{f}(x) = \hat{f}(x) - \widehat{Bias}(\hat{f}(x)),
\end{equation*}
which can be written in the following form
\begin{equation}
\label{rbias}
\tilde{f}(x) = \hat{f}(x) - \frac{h^2}{2}\hat{f}''(x)\int u^2K(u)du = \frac{1}{nh}\sum_{i=1}^n K\Big(\frac{x-X_i}{h}\Big) - \frac{h}{2n}\int u^2K(u)du\sum_{i=1}^nK''\Big(\frac{x-X_i}{h}\Big).
\end{equation}

\bigskip
\noindent 
Further, in this paper, we assume that the hypotheses  $\displaystyle \mathbf{H_1}, \mathbf{H_2}, \mathbf{H_3}$ \, hold for both $K$ and $K''$,  they  are Riemann integrable and twice continuously differentiable, and that  $f$ is bounded with  support included in  $\mathbb{R}_+$. We denote by  $x_0$  the infinimum of this support.
Let's substite in  \, \eqref{cc1} \, $f$ \, by \, $\tilde{f}$. We obtain
\begin{equation} 
\label{cc3}
\tilde{J}_n(z,\alpha)=\int_0^z\Big(\frac{z-x}{z}\Big)^{\alpha}(nh)^{-1}\sum_{i=1}^n K\Big(\frac{x-X_i}{h}\Big)\, - \frac{h}{2n}\int u^2K(u)du\int_0^z\Big(\frac{z-x}{z}\Big)^{\alpha}\sum_{i=1}^nK''\Big(\frac{x-X_i}{h}\Big).
\end{equation}    
Let $\displaystyle \Delta_{hi}=[hi,h(i+1)[, 0\leq i < [\frac{z}{h}]$ be a partition of $[0, z]$ and using the Riemann sum definition of the integral, we establish  that it correspond to the integral $\displaystyle \tilde{J}_n(z,\alpha)$  of the following sum:
\begin{equation} 
\label{cc4}\tilde{P}_n(z,\alpha)=\frac{1}{n}\sum_{j=1}^n\sum_{i=0}^{[z/h]}\Big(\frac{z-ih}{z}\Big)^{\alpha}\Bigg[K\Big(\frac{ih-X_j}{h}\Big)\, - \frac{h^2}{2}\int u^2K(u)duK''\Big(\frac{ih-X_j}{h}\Big)\Bigg]
 + \mathcal{V}_{n,b}(z)
\end{equation}
where  $\displaystyle \Big[\frac{\cdot}{h}\Big]$  design the integer part of  $\displaystyle \frac{\cdot}{h}$ \, and \, $\displaystyle  \mathcal{V}_{n,b}(z) \to 0$ \, in probability as  \, $n \to +\infty$ (for more details, one can refer to Section \ref{const}). Finaly, we propose as  estimator of  FGT poverty index, the following new estimator :
\begin{equation} 
\begin{split}
\label{cc5} P_{n,b}(z,\alpha)  = \frac{1}{n}\sum_{j=1}^n\sum_{i=0}^{[z/h]}\Big(\frac{z-ih}{z}\Big)^{\alpha}\Bigg[K\Big(\frac{ih-X_j}{h}\Big)\, 
 - \frac{h^2}{2}\int u^2K(u)duK''\Big(\frac{ih-X_j}{h}\Big)\Bigg] ,
\end{split}
\end{equation}
$ \text{for} \quad \alpha=0 \quad \text{or} \quad \alpha \geq 1$.

\bigskip
\noindent Additional hypotheses are made about both  $K$ and $K''$, that is:
\begin{enumerate}
\item[$(\mathbf{H_4})$] $K$ \, is of bounded variation function $\displaystyle V_{-\infty}^{u}K$ on $\mathbb{R}$ and let $V(\mathbb{R})$ be its total variation.
\item[$(\mathbf{H_5})$] $\int_{\mathbb{R}}|uK(u)|du < +\infty$ and $\int u^2K(u)du < +\infty$.


\item[$(\mathbf{H_6})$]  There exists a nonincreasing function  $\lambda$  such that $\displaystyle \lambda (\frac{u}{h}) = O(h)$ \, on bounded intervals, 
\[\forall (x, y)\in \mathbb{R}^2, |K(x)-K(y)| \leq \lambda|x-y| \; \text{and} \; \lambda(u) \longrightarrow 0 \; \text{when} \; u \rightarrow 0, \text{and} \; u\geq 0.\]
\end{enumerate}
\section{Main results}
\noindent Our main resultats are  relative to the following additional about the density function  \, $f$:

$\mathbf{C_1} : \quad f$ \, is uniformly continuous.

$\mathbf{C_2} : \quad f$ \, is twice almost everywhere differentiable \, $\displaystyle f', f'' \in L_1(\mathbb{R}) $.
\subsection{Uniform almost sure consistency and behavior of the bias}
\begin{theo}\rm
\label{theo211}
 Assume that the hypotheses  $\mathbf{H_4}$  and $\mathbf{C_1}$ hold. Then for all  $M>0$, estimator $\displaystyle P_{n,b}(z, \alpha)$ converges uniformly almost surely on  $[0, M]$  to   $\displaystyle P(z, \alpha)$ as $\displaystyle n \rightarrow +\infty \quad  i.e.$
 \[P\Big(\lim_{n \to +\infty}\sup_{z \in [0, M]}|P_{n, b}(z, \alpha)-P(z, \alpha)|=0\Big)=1,\]
provided $\displaystyle nh^2(\log\log n)^{-1} \to +\infty$  as  $\displaystyle n \to +\infty$.
\end{theo}
\begin{theo}\rm
\label{theo221}
  Assume that the hypotheses   $\mathbf{H_4}, \mathbf{H_5}$  and  $\mathbf{C_2}$ hold. Then for all  $M>0$, the estimator $\displaystyle P_{n, b}(z, \alpha)$ converges uniformly almost surely on  $[0, M]$  to  $\displaystyle P(z, \alpha)$ as $\displaystyle n \rightarrow +\infty \quad i.e.$
 \[P\Big(\lim_{n \to +\infty}\sup_{z\in [0, M]}|P_{n, b}(z, \alpha)-P(z, \alpha)|=0\Big)=1,\]
 provided that  $\displaystyle nh^2(\log\log n)^{-1} \to +\infty$ as $n \to +\infty$.
\end{theo}  
For the demonstration of the theorems, we use the Theorem 2 of Kiefer \cite{Kiefer1961} and the following lemmas showing that  $P_{n, b}(z, \alpha)$ is uniformly asymptotic unbiased on all bounded interval.
\begin{lem}\rm
\label{lem211}
 If the hypothese  $\mathbf{C_1}$ holds , then $\forall \,  M>0$, we have
 \[\lim_{n \to +\infty}\sup_{z\in [0, M]}|\mathbb{E}(P_{n,b}(z, \alpha))-P(z, \alpha)|= 0, .\]
\end{lem}
\begin{lem}\rm
\label{lem221}
 If the hypotheses $\mathbf{H_5} \; \text{and} \; \mathbf{C_2}$ hold, then:
 \begin{equation}
 \begin{split}
 &\sup_{z\in \mathbb{R}}|\mathbb{E}(P_{n,b}(z, \alpha))-P(z, \alpha)| \leq h\Bigg(\Big(\int_{\mathbb{R}}|f^{'}(x)|\, dx\Big)\Big(\int_{\mathbb{R}}(|u|+1)|K(u)|\, du\Big)\\
 &+2(\alpha D+Ah)\int_{-\infty}^{+\infty}|K(u)|\, du\Bigg) + h^3\int u^2K(u)du,
 \end{split}
 \end{equation}
 where
 \[D=\sup_{x\in\mathbb{R}}F(x) \quad \text{and} \quad A=\sup_{x\in\mathbb{R}}f(x).\]
 \end{lem}
\begin{rem}\rm
  If  $K$ satisfies the hypothese  $\mathbf{H_5}$, then by using  $\mathbf{H_1}$, the kernel

   $\displaystyle \widehat{K}=\frac{K^2}{\int_{\mathbb{R}}K^2(y)\, dy}$  also satisfy it. 
 \end{rem}
From the two previous lemmas, we get  the following corollaries:
 \begin{cor}\rm
 \label{cor211}
 Under the assumptions of Lemma \ref{lem211}, we have uniformly on $[0, M]$ (resp $\mathbb{R}$)
\begin{equation*}
\begin{split}
\lim_{n\to+\infty}\mathbb{E}\Bigg(\sum_{i=1}^{[\frac{z}{h}]}\Big(1-\frac{ih}{z}\Big)^{2\alpha}\Bigg[K\Big(\frac{ih-X_j}{h}\Big) - \frac{1}{2}\int u^2K(u)du K''\Big(\frac{ih-X_j}{h}\Big)\Bigg]^2 \Bigg) = \Big(\int_{\mathbb{R}}K(y)\, dy\Big)P(z, 2\alpha).
\end{split}
\end{equation*}
 \end{cor}
 \begin{cor}\rm
 \label{cor221}
 If the assumptions of Theorem \ref{theo221} hold and if $\displaystyle h= O(n^{-1}\log\log n)^{1/4}$, then for all $M>0$, we have almost surely: 
  \[\sup_{z\in [0, M]}|\mathbb{E}(P_{n,b}(z, \alpha))-P(z, \alpha)|= O(n^{-1}\log\log n)^{3/4}.\]
 \end{cor}
 \subsection{Uniforme  mean square consistency}
 \begin{theo}\rm
 \label{theo231}
 If   $\mathbf{H_6}$  and  $\mathbf{C_1}$ hold. Then:
 \begin{enumerate}    
\item $\displaystyle \lim_{n \to +\infty}n\mathbb{V}ar(P_{n,b}(z, \alpha))=\Big(\int_ {\mathbb{R}}K^2(y)\, dy\Big)P(z, 2\alpha)-\Big(P(z, \alpha)\Big)^2$.
 \item $\text{For all} \quad M>0$,
  \[\lim_{n \to +\infty}\sup_{z\in [0, M]}\mathbb{E}\Big(P_{n, b}(z, \alpha)-P(z, \alpha)\Big)^2=0.\]
 \end{enumerate}
 \end{theo}
 \begin{theo}\rm
 \label{theo241}
 Assume that  $\mathbf{H_6}$  and  $\mathbf{C_2}$ hold. Then:
 \begin{enumerate}    
\item $\displaystyle  \lim_{n \to +\infty}n\mathbb{V}ar(P_{n, b}(z, \alpha))=\Big(\int_ {\mathbb{R}^2}K^2(y)\, dy\Big)P(z, 2\alpha)-\Big(P(z, \alpha)\Big)^2$.
\item $ \text{Moreover, if } \quad \mathbf{H_5} \quad \text{holds, we have \, for all} \quad M>0$,
 \[\lim_{n \to +\infty}\sup_{z\in [0, M]}\mathbb{E}\Big(P_{n, b}(z, \alpha)-P(z, \alpha)\Big)^2=0.\]
\end{enumerate}    
\end{theo}
\noindent For the proof of this theorem, we assume that the hypothese $\mathbf{C_1}$  or  $\mathbf{C_2}$ holds and before that, we prove the  Theorem \ref{theo251} below by using the following lemma.
 \begin{lem}\rm
 \label{lem231} Let $0 \leq \theta_i \leq 1, i=1,2$. Then for all  $x, y$  and $x \neq y$ we have
  \begin{equation*}
  \begin{split}
  \lim_{n \to +\infty}\Big((h)^{-2}\int_{-\infty}^{+\infty} &\Bigg(\Bigg|\Bigg[K\big(\frac{u-x+\theta_1}{h}\big) - \frac{1}{2}\int s^2K(s)ds K''\big(\frac{u-x+\theta_1}{h}\big)\Bigg]\\
  &\times\Bigg[K\big(\frac{u-y+\theta_2}{h}\big) - \frac{1}{2}\int s^2K(s)ds K''\big(\frac{u-y+\theta_2}{h}\big)\Bigg]\Bigg|
  \Bigg)f(u)\, du\Big)=0.
  \end{split}
  \end{equation*}
 \end{lem}
\begin{theo}\rm
\label{theo251}
Assume that the hypothese  $\mathbf{H_6}$ holds. Then for all $M>0$,
\begin{equation*}
  \begin{split}
\lim_{n \to +\infty}\sup_{z\in [0, M]}\sum_{0\leq i\neq j\leq [\frac{z}{h}]}\Big(1-\frac{ih}{z}\Big)^{\alpha}\Big(1-\frac{jh}{z}\Big)^{\alpha}&\int_{\mathbb{R}}\Big[K\Big(\frac{u - ih}{h}\Big) - \frac{h^2}{2}\int s^2K(s) ds K''\Big(\frac{u - ih}{h}\Big)\Big]\\
& \times\Big[K\Big(\frac{u - jh}{h}\Big) - \frac{h^2}{2}s^2K(s)dsK''\Big(\frac{u - jh}{h}\Big)\Big]f(u)\, du=0.
\end{split}
\end{equation*}
\end{theo}
\begin{rem}\rm
The estimator $P_{n, -b}(z,\alpha)$  has asymptotic efficiency with respect to $\widehat{P}_n(z,\alpha)$,  $$e(z,\alpha)=\bigg((\int_{\mathbf{R}}K^2(y)\,dy) P(z,2\alpha)-(P(z,\alpha))^2\bigg)/P(z,2\alpha)-(P(z,\alpha))^2.$$
 The integral $\int_{\mathbf{R}}K^2(y)\,dy$  is strictly less than 1 for the conventional kernels \cite{Parzen1962} p.1068. Then we have in this case $e(z,\alpha)<1$. In theorem \ref{theo241}, the speeed of convergence in mean square is of the order of $O(\frac{1}{n^2})$ if $h$ is of the order of $O(\frac{1}{\sqrt{n}})$.  
\end{rem}
\section{Simulation Study}
\noindent In this section, we make a simulation study giving the mean square error and empirical variance of 50 samples of size $n=1000$ of the three estimators that we compared. Our bias reduction kernel estimator and the classical one are evaluated by a Gaussian kernel checking assumptions $H_i , i = 1, ..., 6$, taking $h = (n log n)^{-1/2}$.
For a Pareto distribution type on $[0,1]$ with parameters  $x_0= 0.02$   and $ \beta=0.2$
, we calculated the mean square error $mse_1$ of $\displaystyle (P_{n,b, 1}(z,\alpha), ...,  P_{n,b, 50}(z, \alpha))$,  $mse_2$ of $\displaystyle (P_{n, 1}(z,\alpha), ...,  P_{n, 50}(z, \alpha))$ and $mse_3$ of $\displaystyle (P_{n,1}^\lambda(z, \alpha), ..., P_{n, 50}^\lambda(z, \alpha))$ and the respective empirical variances $\sigma_1^2$, $\sigma_2^2$ and  $\sigma_3^2$  for different values of $(z,\alpha)$ by the following statistics:
$\overline{P_{n,b}(z,\alpha)}=\frac{1}{50}\sum_{i=1}^{50}P_{n,b,i}(z,\alpha)$,
$mse_1=\frac{1}{50}\sum_{i=1}^{50}(P_{n,b,i}(z,\alpha)-P(z,\alpha))^2$ 
and
$\sigma_1^2=\frac{1}{50}\sum_{i=1}^{50}\left(P_{n,b,i}(z,\alpha)-\overline{P_{n,b}(z,\alpha)} 
\right)^2.
$ Similarly, $\overline{P_n(z,\alpha)}$ and $\overline{P_n^\lambda(z,\alpha)}$, $(mse_2,\sigma_2^2)$  and $(mse_3, \sigma_3^2)$ are respectively calculated  for the estimator $P_n(z,\alpha)$ and $P_n^\lambda(z,\alpha)$. The results are ????????
$$
\begin{tabular}{*{9}{c}}
\hline
\hline
\multicolumn{8}{c}{}\\
z&0.1& 0.2& 0.3 & 0.4 & 0.5& 0.6 & 0.7 \\
\multicolumn{8}{c}{}\\
\hline
\hline
$\alpha=0$& \multicolumn{7}{c}{ }\\
$mse_1$&0.005562&0.155562&0.188334&0.261614&0.212470&0.244210&0.251033\\
$mse_2$&0.125552&0.165542&0.198308&0.271581&0.218349&0.254173&0.260993\\
$mse_3$&0.011194&0.160774&0.206454&0.277992&0.247950&0.264731&0.281827\\
$\sigma_1^2$&0.0002778&0.0005931&0.0006942&0.0008862&0.00077324&0.0008681&0.0007603\\
$\sigma_2^2$&0.0002779&0.0005932&0.0006943&0.0008864&0.00077329&0.0008683&0.0007605\\
$\sigma_3^2$& 0.0051733&0.0096222&0.0097822&0.0099378&0.01013778&0.0103822&0.01073778\\
\hline
& \multicolumn{7}{c}{ }\\
\hline
$\alpha=1$& \multicolumn{7}{c}{ }\\
$mse_1$&2.551305&2.335935&2.213879&2.133226&2.07242&2.024278&2.073725\\
$mse_2$&2.651286&2.435894&2.313825&2.233163&2.17235&2.124203&2.083646\\
$mse_3$&2.373474&2.413784&2.291902&2.217399&2.161699&2.116292&2.080483\\
$\sigma_1^2$&5.6511e-05&0.00021540&0.0003155&0.00039351&0.0004510&0.0005005&0.0005335\\
$\sigma_2^2$&5.6521e-05&0.00021544&0.0003156&0.00039359&0.0004511&0.0005006&0.0005336\\
$\sigma_3^2$ &0.000996&0.003717993&0.005432&0.0061841&0.006786&0.007338&0.007711\\
\hline
& \multicolumn{7}{c}{ }\\
\hline
$\alpha=2$& \multicolumn{7}{c}{ }\\
$mse_1$&2.573554&2.176867&2.031041&2.028172&1.871536&1.820177&1.778642\\
$mse_2$&2.573544&2.276842&2.131005&2.038128&1.971485&1.920121&1.878582\\
$mse_3$&2.587938&2.266043&2.113916&2.020992&1.956386&1.907174&1.867879\\
$\sigma_1^2$&1.737e-05&9.920e-05&0.0001805&0.0002456&0.0003002&0.0003469&0.0003871\\
$\sigma_2^2$& 1.738e-05&9.921e-05&0.0001806&0.0002457&0.0003003&0.000347&0.0003872\\
$\sigma_3^2$ &0.0002843&0.001747&0.003194&0.004243&0.0049853&0.00557&0.006058097\\
\hline
\hline
\end{tabular}
$$
The studies cases $\displaystyle {P}(z, 0), {P}(z_, 1), {P}(z, 2)$ are commonly and respectively called the poverty rate or headcount ratio, the depth of poverty or poverty gap index and the severity of poverty \cite{Foster&Greer&Thorbecke1984}. 
%
\noindent A comparison of simulation results shows that for small samples, each point $z$, our bias reduction kernel estimator provides a much lower error and variance for the three values of $\alpha$ considered. Thus, we can conclude that our estimator is recommended.
\section{Details of the Proofs}
\subsection{Construction of the estimator}
\label{const}
For \, $z > 0$ \, and
$\Delta_{h,i} = [hi, h(i+1)[ \qquad i=0, \cdot \cdot \cdot,\big[\frac{z}{h}\big]$.
We have the following Riemann sum over the interval \, $\left[0, z\right]$:
\begin{equation*}
\begin{split}
S_n(z,\alpha)=\frac{1}{n}&\sum_{j=1}^n\sum_{i=0}^{[\frac{z}{h}]-1}\Big(1-\frac{ih}{z}\Big)^{\alpha}\Bigg[K\Big(\frac{ih-X_j}{h}\Big) - \frac{h^2}{2}\int u^2K(u)duK''\Big(\frac{ih-X_j}{h}\Big)\Bigg]\\
&+\Big(z-h[\frac{z}{h}]\Big)\sum_{j=1}^n\Bigg(1-\frac{h[\frac{z}{h}]}{z}\Bigg)^{\alpha}\frac{1}{n}\frac{1}{h}\Bigg[K\Bigg(\frac{[\frac{z}{h}]h-X_j}{h}\Bigg) - \frac{h^2}{2}\int u^2K(u)duK''\Big(\frac{[\frac{z}{h}]h-X_j}{h}\Big)\Bigg],
\end{split}
\end{equation*}
corresponding to the integral
\[J_n(z,\alpha)=\int_0^{z}\Big(\frac{z-x}{z}\Big)^{\alpha}\frac{1}{n}\sum_{j=1}^n\frac{1}{h}\Bigg[K\Big(\frac{x-X_j}{h}\Big) -  \frac{h^2}{2}\int u^2K(u)duK''\Big(\frac{x- X_j}{h}\Big)\Bigg]\, dx.\]
The sum $S_n$ can be rewritten as
\begin{align*}
S_n=\frac{1}{n}\sum_{j=1}^n\sum_{i=0}^{[\frac{z}{h}]}\Big(1-\frac{ih}{z}\Big)^{\alpha}\Bigg[K\Big(\frac{ih-X_j}{h}\Big) - \frac{h^2}{2}\int u^2K(u)duK''\Big(\frac{ih-X_j}{h}\Big)\Bigg] + \mathcal{V}_{n,b}(z),
\end{align*}
with
\begin{align*}
\mathcal{V}_{n,b}(z) =\frac{1}{n}\sum_{j=1}^n\frac{\Big(z-h[\frac{z}{h}]\Big)-h}{h}\Big(1-\frac{h[\frac{z}{h}]}{z}\Big)^{\alpha}\Bigg[K\Big(\frac{[\frac{z}{h}]h-X_j}{h}\Big) - \frac{h^2}{2}\int u^2K(u)duK''\Big(\frac{[\frac{z}{h}]h-X_j}{h}\Big)\Bigg].
\end{align*}
Now, we have to show that
$\mathcal{V}_{n, b} \longrightarrow 0$ in probability as $n \rightarrow \infty$.
Since $[\frac{z}{h}]\leq \frac{z}{h}  <[\frac{z}{h}]+1 $, we get for $\alpha\geq 0$ :
\begin{equation}
\label{v11}
|\mathcal{V}_{n,b}(z)|\leq 
\frac{1}{n}\sum_{j=1}^n\Bigg|K\left(\frac{[\frac{z}{h}]h-X_j}{h}\right) - \frac{h^2}{2}\int u^2K(u)duK''\left(\frac{[\frac{z}{h}]h-X_j}{h}\right)\Bigg|.
\end{equation}
Note that for all $x\in \mathbb{R}$,  $hx+h[\frac{z}{h}]\longrightarrow z$ as $n\longrightarrow \infty$. Thus, for $\alpha = 0$, using the continuity of $f$ and the fact that $K,K^{''} \in L_1(\mathbb{R})$ , we have for $n$ large enough:
\begin{equation} 
\begin{split}
\Big|\mathbb{E}\Big[\frac{1}{n}\sum_{k=1}^n\Bigg\{K\Big(\frac{[\frac{z}{h}]h-X_j}{h}\Big) -  \frac{h^2}{2}\int s^2K(s)ds&K''\Big(\frac{[\frac{z}{h}]h-X_j}{h}\Big)\Bigg\}\Big]\Big|\\
 & = \Big|h\int \Bigg[K(u)- \frac{h^2}{2}\int s^2K(s)dsK''(u)\Bigg]f(h[\frac{z}{h}]-uh)du\Big|\\
& \leq  hf(z)(1+o(1))\int |K(u)|du.
\end{split}
\end{equation}
Hence, $\mathbb{E}|\mathcal{V}_{n,b}(z)|=O(h)$ as $n\longrightarrow \infty$. The estimator holds by using the markov's inequality.
%
%
%
%
%
%
\subsection{Proofs of main results}
\subsubsection*{Proof of lemma \ref{lem211}}
\noindent Let define $\displaystyle \bar{\Delta}_{h,i}=\Delta_{h,i}\cap[0, M];$ \, and \, $\chi_B$ \, the indicator function of the set \, $B$.
Put $$\varphi_h(u)=\frac{h^2}{2}\int s^2K(s)dsK''(u).$$
By a change of variables, we have
\begin{eqnarray*}
\mathbb{E}(P_n(z, \alpha))-P(z,\alpha)&=& \sum_{i=1}^{[\frac{z}{h}]}h\Big(1-\frac{ih}{z}\Big)^{\alpha}\int (K(u) -\varphi_h(u))f(ih-uh)\, du-P(z,\alpha)\\
&=&\int_0^{z}\sum_{i=1}^{[\frac{z}{h}]}\chi_{\bar{\Delta}_{h,i}}(x)\Big(1-\frac{ih}{z}\Big)^{\alpha}\int(K(u) - \varphi_h(u))f(ih-uh)\, dudx-P(z,\alpha)\\
&\;&\quad +\; (h([\frac{z}{h}]+1)-z)\Big(1-\frac{h[\frac{z}{h}]}{z}\Big)^{\alpha}\int(K(u) - \varphi_h(u)\Bigg)f(h[\frac{z}{h}]-uh)\, du\\
&:=& A_{n,1}+A_{n,2}.
\end{eqnarray*}
We first study the term $A_{2,n}$. Let $z \in [0, M]$, Since 
$\displaystyle |h([\frac{z}{h}]+1)-z| \leq h$ and $0<h[\frac{z}{h}]\leq z$, we get for $n$ large enough:
\begin{eqnarray*}
\label{31}
\sup_{z\in[0,M]}|A_{n,2}|&:=&\sup_{z\in[0,M]}\left|(h([\frac{z}{h}]+1)-z)\Big(1-\frac{h[\frac{z}{h}]}{z}\Big)^{\alpha}\int(K(u)-\varphi_h(u))f(h[\frac{z}{h}]-uh  )\, du \right|\\
&\leq& h(\sup_{z\in [0,M]}f(z))(1+o(1))\int |K(u)| du.
\end{eqnarray*}
Now, we are going to study the term $A_{1,n}$. Under $\mathbf{H_2}$, one can rewrite $P(z, \alpha)$ as:
$$
P(z, \alpha)= \int_0^{z}\Big(1-\frac{x}{z}\Big)^{\alpha}\int (K(u) -\varphi_h(u))\, du f(x)dx+\int_0^{z}\Big(1-\frac{x}{z}\Big)^{\alpha}\int \varphi_h(u)\, duf(x)dx.
$$
Let $\displaystyle x \in \bar{\Delta}_{h,i}, \; i=0,...,[z/h]$. We have 
\begin{eqnarray} 
\label{51}
\nonumber &\;& \left|\Big(1-\frac{ih}{z}\Big)^{\alpha}\int (K(u) -\varphi_h(u)) f(ih-uh)\, du-\Big(1-\frac{x}{z}\Big)^{\alpha}\int (K(u) -\varphi_h(u))\, duf(x)\right|\\
 \nonumber  &\;&\quad\quad \quad\quad\quad=\left|\int\Big[\Big(1-\frac{ih}{z}\Big)^{\alpha}f(ih-uh)-\Big(1-\frac{x}{z}\Big)^\alpha
  f(x)\Big](K(u) -\varphi_h(u))\, du\right |\\
  &\;&\quad\quad \quad\quad\quad \leq \quad \quad \int \left|\Big(1-\frac{ih}{z}\Big)^{\alpha}- \Big(1-\frac{x}{z}\Big)^{\alpha}\right|
  \times \left|K(u) -\varphi_h(u)\right||f(x)\, du\\
 \nonumber &\;&\quad\quad \quad\quad\quad\quad + \quad\quad\int \Big(1-\frac{ih}{z}\Big)^{\alpha}\left|f(ih-uh)-f(x)\right|\times\left|K(u)-\varphi(u)\right|duf(x).   
 \end{eqnarray} 
For $\alpha=0$, we have $ \Big|\Big(1-\frac{ih}{z}\Big)^{\alpha}-\Big(1-\frac{x}{z}\Big)^{\alpha}
   \Big|=0$. Next, for $\alpha\geq 1$, using 
the first order Taylor formula on $[ih,x]$ to the function 
$g(t)= \Big(1-\frac{t}{z}\Big)^{\alpha}$, there exists $c_i \in ]hi, x[$, such that  
$g(ih)-g(x)=g^\prime(c_i)(ih-x)$. That is  
\begin{equation*}
   \begin{split}   
   \Big|\Big(1-\frac{ih}{z}\Big)^{\alpha}-\Big(1-\frac{x}{z}\Big)^{\alpha}
   \Big|& = \frac{\alpha }{z}\Big|1-\frac{c_i}{z}\Big|^{\alpha-1}\Big|ih - x\Big|\leq \frac{\alpha h}{z}.
\end{split}
\end{equation*}
Denote by  $\displaystyle I_1^{i}(x)$ ( respectively $\displaystyle I_2^{i}(x)$),  the first term (respectively the second term ) of the right hand-side of the inequality \eqref{51}. 
For simplify the notations, let 
\[
I_1(x) = \sum_{i=1}^{[\frac{z}{h}]}\chi_{\bar{\Delta}_{h,i}}(x)I_1^{i}(x)
\quad\text{and}\quad
I_2(x) = \sum_{i=1}^{[\frac{z}{h}]}\chi_{\bar{\Delta}_{h,i}}(x)I_2^{i}(x).
\]
We have 
$$|A_{1,n}|\leq\int_0^z I_1(x)dx+\int_0^z I_2(x)dx+
\int_0^{z}\Big(1-\frac{x}{z}\Big)^{\alpha}\int |\varphi_h(u)|\, duf(x)dx.
$$
We remark by definition that $\int |\varphi_h(u)|\, du=O(h^2)$. This leads to 
$$
|A_{1,n}|\leq\int_0^z I_1(x)dx+\int_0^z I_2(x)dx+O(h^2).
$$
Now, we are going to study the terms $\int_0^z I_1(x)dx$ and $\int_0^z I_2(x)dx$.
First, we have  
\begin{equation}
\label{61}
\int_0^{z}I_1(x)\, dx \leq \frac{\alpha h}{z}\int_0^{z}(\int f(x)|K(u) - \varphi_h(u)|\, du)\, dx
\leq\frac{\alpha h}{z}F(z)\int \left(|K(u)|+|\varphi_h(u)|\right)\, du.
\end{equation} 
Remarking also that $\int |\varphi_h(u)|\, du=O(h^2)$, we get by assumptions
$$\int_0^{z}I_1(x)\, dx=O(h).$$
Next, note that
\begin{equation}
\begin{split}
\label{71}
I_2^{i}(x)&\leq \Bigg(\int|f(ih-uh)-f(ih)|\times|K(u) -\varphi_h(u)|\, du\\
& \qquad +\int |f(ih)-f(x)|\times|
K(u) - \varphi_h(u)|\, du\Bigg) .  
\end{split}
\end{equation}
Let $\varepsilon>0$, since $f$ is uniformly continuous, there exists $\displaystyle \eta_0=\eta_0(z)>0$  such that $\displaystyle |ih-x|\leq \eta_0$  hence if  $\displaystyle h\leq \eta_0$, we have  $\displaystyle |f(ih)-f(x)|<\frac{\varepsilon}{M}$. Therefore, 
\begin{equation*}
\begin{split}
\int_0^{z}I_2(x)\, dx & \leq \sum_{i=1}^{[\frac{z}{h}]}h\int f(ih-uh)-f(ih)|\times|K(u) - \varphi_h(u)|\,du + \varepsilon\int |K(u) - \varphi_h(u)|\, du.  
\end{split}
\end{equation*}   
By the uniform continuity of $f$  we have
 \[\exists \eta_1=\eta_1(z)>0, |uh|<\eta_1 \Rightarrow |f(uh-ih)-f(ih)|< \frac{\varepsilon}{M}.\]
Hence,
\begin{equation*}
\begin{split}
&\sum_{i=1}^{[\frac{z}{h}]}h\int |f(ih-uh)-f(ih)|\times|K(u) - \varphi_h(u)|\, du\\
&\leq \sum_{i=1}^{[\frac{z}{h}]}h\int_{|uh|< \eta_1}\frac{\varepsilon}{M}|K(u) - \varphi_h(u)|\, du+ \sum_{i=1}^{[\frac{z}{h}]}h\int_{|uh|\geq\eta_1}|f(ih-uh)-f(ih)|\times|K(u) -\varphi_h(u)|\, du \\
&\leq \varepsilon\int |K(u) - \varphi_h(u)|
\,du+\sum_{i=1}^{[\frac{z}{h}]}h\int_{|uh|\geq \eta_1}\times|f(ih-uh)|\times|K(u) -\varphi_h(u)|\,du\\
& \qquad \qquad   + \sum_{i=1}^{[\frac{z}{h}]}h\int_{|uh|\geq \eta_1}f(ih)\times|K(u) -\varphi_h(u)|\, du. 
\end{split}
\end{equation*}  
Since $f$ is continuous, it is Riemann-integrable and by remarking that 
$(h([\frac{z}{h}]+1)-z)f(h[\frac{z}{h}]) \to 0$, $n \to +\infty$, we have the sum
$\sum_{i=1}^{[\frac{z}{h}]}hf(ih) \to \int_0^{z}f(x)\, dx=F(z) $
is bounded. Let  $A$  be the latter.\\
By the change of variables $\displaystyle v=uh$, we get
\begin{equation*}
\begin{split}
&\sum_{i=1}^{[\frac{z}{h}]}h\int_{|uh|\geq \eta_1}f(ih-uh)|K(u) - \varphi_h(u)|\, du
+\sum_{i=1}^{[\frac{z}{h}]}h\int_{|uh|\geq \eta_1}f(ih)|K(u)-\varphi_h(u)|\, du \\
&\leq\sum_{i=1}^{[\frac{z}{h}]}h\int_{|v|\geq \eta_1}f(ih-v)\frac{1}{h}|K(\frac{v}{h}) - \varphi(\frac{v}{h}) |\, dv +A\int_{|uh|\geq \eta_1}|K(u) -\varphi_h(u)|\, du. 
\end{split}
\end{equation*}
Under $(\mathbf{H_3})$, there exists  $C>0$  fixed, such that  
$\displaystyle \Big|\frac{v}{h}\Big| \geq C$  we have 
 \[\Big|\frac{v}{h}\Big|\Big|K(\frac{v}{h}) -\varphi_h(\frac{v}{h})\Big|\leq \frac{\eta_1\varepsilon}{M}.\]
Let  $\displaystyle \eta = \inf(\eta_1, Ch)=Ch, h$  being small enough,
then
\begin{equation*}
\begin{split}
\frac{1}{\eta_1}\int_{|v|\geq \eta_1}|f(ih-v)|\frac{v}{h}|K(\frac{v}{h}) -\varphi_h(\frac{v}{h})|\, dv&\leq \frac{\varepsilon}{M}\int_{|v|\geq \eta_1}f(ih-v)\, dv\\
&\leq \frac{\varepsilon}{M}\int_{\mathbb{R}}f(x)\, dx = \frac{\varepsilon}{M} .
\end{split}
\end{equation*}
Hence,
\[\sum_{i=1}^{[\frac{z}{h}]}h\int_{|v|\geq \eta_1}f(ih-uh)\frac{1}{h}|K(\frac{v}{h}) -\varphi_h(\frac{v}{h})|\, dv\leq \frac{z\varepsilon}{M}\leq \varepsilon. \]
Since $\displaystyle A\int_{|uh|\geq \eta_1}|K(u) -\varphi_h(u)|\, du \to 0, n \to +\infty$, we have together with \eqref{61}
\begin{equation*}
\begin{split}
\lim_{n \to +\infty}\sup_{z\in [0, M]}|\mathbb{E}(P_n(z, \alpha))-P(z, \alpha)|\leq 2\varepsilon\int|K(u)|, du+\varepsilon.   
\end{split}
\end{equation*}
The proof of the lemma is complete.
\subsection*{Proof of lemma \ref{lem221}}
 \noindent For  $x \in \bar{\Delta}_{h,i}$, $i=1, \cdot \cdot \cdot, [\frac{z}{h}]$. We have
\begin{equation*}
\begin{split}
\int \Bigg(1-\frac{ih}{z}\Bigg)^{\alpha}&|f(ih-uh)-f(x)|\times|K(u) - \varphi_h(u)|\, du\\
&\leq\Bigg(\int \int_{x}^{x+h(|u|+1)}|f^{'}(t)|\times|K(u) - \varphi_h(u)|\, dtdu\Bigg)  
\end{split}
\end{equation*} 
Hence, by using again the expression of  $I_2(x)$ defined in the proof of Lemma \ref{lem211}, we have
\begin{equation}
\begin{split} 
 \label{81}
 \int_0^{z}I_2(x)\, dx\leq \int_0^{z}\Big(\int \int_{x}^{x+h(|u|+1)}|f^{'}(t)||K(u) -\varphi_h(u)|\, dtdu\Big)\, dx.
\end{split}
\end{equation}
By the change of variable with $\displaystyle t=x+h(|u|+1)v$ and by using Fubini's theorem,weget
\begin{equation*}
\begin{split}
\int_0^{z}I_2(x)\, dx&\leq h\int(|u|+1)|K(u) -\varphi_h(u)|du\\
& \quad \times \Big(\int |f^{'}(x+h(|u|+1)v)|\times|K(u) - \varphi_h(u)|\, du\Big)\, dx\, du\int_0^1\, dv\\
&=\qquad  \int_{\mathbb{R}}|f^{'}(x)|\, dx.
\end{split}
\end{equation*} 
 This inaquality together with \eqref{31} and \eqref{61} lead to completion of the proof.
\subsection*{Proof of theorem \ref{theo211} and theorem \ref{theo221}}
 \begin{preuve}
\noindent Let  $\displaystyle \widehat{F}_n$  be the empirical distribution of the sample $\displaystyle (X_1, X_2 , \cdot, \cdot, \cdot, X_n)$  defined by
 \[F_n(l)=n^{-1}\sum_{i=1}^n(\chi_{X_i<l})\] 
where  $\displaystyle \chi_A$ stands for the indicator function  $A$. We can write
\[P_n(z, \alpha)=\int_{\mathbb{R}}\sum_{i=1}^{[\frac{z}{h}]}\Big(1-\frac{ih}{z}\Big)^{\alpha}\Bigg[K\Big(\frac{l-ih}{h}\Big) - \frac{h^2}{2}\int s^2K(s)dsK''\Big(\frac{l-ih}{h}\Big)\Bigg]\, d\hat{F}_n(l)\]
and
\[\mathbb{E}(P_n(z, \alpha))=\sum_{i=1}^{[\frac{z}{h}]}\int_{\mathbb{R}}\Big(1-\frac{ih}{z}\Big)^{\alpha}\Bigg[K\Big(\frac{l-ih}{h}\Big)  - \frac{h^2}{2}\int s^2K(s)dsK''\Big(\frac{l-ih}{h}\Big)\Bigg]\, dF(l).\]
 We have  
 \begin{equation*}
\begin{split}
 |P_n(z, \alpha)-\mathbb{E}(P_n(z, \alpha))|=\Bigg|\sum_{i=1}^{[\frac{z}{h}]}\int_{\mathbb{R}}\Big(1-\frac{ih}{z}\Big)^{\alpha}\Bigg[K\Big(\frac{l-ih}{h}\Big) - \frac{h^2}{2}\int s^2K(s)dsK''\Big(\frac{l-ih}{h}\Big)\Bigg]\, (d\hat{F}_n(l)-dF(l))\Bigg|
\end{split}
\end{equation*}
The integration by parts yields
\begin{equation*}
\begin{split}
|P_n(z, \alpha)-\mathbb{E}(P_n(z, \alpha))| &\leq\sum_{i=1}^{[\frac{z}{h}]}\int_{\mathbb{R}}\Bigg|d\Bigg[K\Big(\frac{l-ih}{h}\Big) -   \frac{h^2}{2}\int s^2K(s)dsK''\Big(\frac{l-ih}{h}\Big)\Bigg]\Bigg|\sup_{l\in \mathbb{R}}|\hat{F}_n(l)-F(l)| \\ 
&\leq\sum_{i=1}^{[\frac{z}{h}]}\int_{\mathbb{R}}\, d\Bigg[K_{-\infty}^{\frac{l-ih}{h}} -  \frac{h^2}{2}\int s^2K(s)dsK{''}_{- \infty}^{\frac{l-ih}{h}}\Bigg]\sup_{l\in\mathbb{R}}|\hat{F}_n(l)-F(l)| \\
&\leq \Big[\frac{z}{h}\Big]V(\mathbb{R})\Big(1 - \frac{h^2}{2}\int s^2K(s)ds\Big) \sup_{l\in\mathbb{R}}|\hat{F}_n(l)-F(l)|     
\end{split}
\end{equation*}
Remarking that 
\begin{equation*}
\begin{split}
|P_n(z, \alpha)-(P(z, \alpha))|&=|P_n(z, \alpha)-\mathbb{E}(P_n(z, \alpha))+\mathbb{E}(P_n(z, \alpha))-P(z, \alpha)|\\
&\leq |P_n(z, \alpha)-\mathbb{E}(P_n(z, \alpha))|+|\mathbb{E}(P_n(z, \alpha))-P(z, \alpha)|
\end{split}
\end{equation*}
By lemma \ref{lem211} we have
\[|\mathbb{E}(P_n(z, \alpha))-P(z, \alpha)| \to 0 \qquad n \to +\infty\]
and  the previous results we have
\[|P_n(z, \alpha)-\mathbb{E}(P_n(z, \alpha))|\to 0 \qquad n \to +\infty\]
\end{preuve}

\subsection{The uniforme mean square consistency}
\subsection*{Proof lemma \ref{lem231}}
\begin{preuve}
We suppose that $\mathbf{C_1}$ holds. Let $\displaystyle \delta>0$.\\
Define
\begin{equation*}
\begin{split}
I_n(x,y)& =(h)^{-2}\int_{-\infty}^{+\infty}\Bigg|\Bigg[K(\frac{u-x+\theta_1h}{h}) - \frac{h^2}{2}\int s^2K(s)dsK''(\frac{u-x+\theta_1h}{h})\Bigg] \\
&\times \Bigg[K(\frac{u-y+\theta_2h}{h})- \frac{h^2}{2}\int s^2K(s)dsK''(\frac{u-y+\theta_2h}{h})\Bigg]|f(u)\, du\\
&=\int_{-\infty}^{+\infty}\Bigg((h)^{-1}\Bigg[K(\frac{v}{h}) - \frac{h^2}{2}\int s^2K(s)dsK''(\frac{v}{h})\Bigg]\Bigg)\\
& \times \Bigg|(h)^{-1}\Bigg[K(\frac{v+x-\theta_1-y+\theta_2h}{h}) - \frac{h^2}{2}\int s^2K(s)dsK''(\frac{v+x-\theta_1-y+\theta_2h}{h})\Bigg|f(x+v-\theta_1h)\, du\\
&=\int_{|v-\theta_1h|\leq \delta}+\int_{|v-\theta_1h|> \delta}
\end{split}
\end{equation*}
 Since  $f(x)$ is continuous, it is bounded on  $\displaystyle I=[x-\delta,x+\delta]$. We assume  $n$ large enough such that $\displaystyle x+v\pm\theta_1h\in I$. Therefore
\begin{equation}
\begin{split} 
\label{91}
\int_{|v-\theta_1h|\leq \delta}&\leq \sup_{|v-\theta_1h|\leq \delta}f(x+v-\theta_1h)\int_{-\frac{\delta}{h}+\theta_1\leq u\leq \frac{\delta}{h}+\theta}|K(u) - \frac{h^2}{2}\int s^2K(s)dsK''(u)|\\
& \times |K(\frac{x-\theta_1h-y+\theta_2h}{h}+u) - \frac{h^2}{2}\int s^2K(s)dsK''(\frac{x-\theta_1h-y+\theta_2h}{h}+u)|(h)^{-1}\, du
\end{split}
\end{equation} 
\begin{equation*}
\begin{split}
&=\sup_{|v-\theta_1h|\leq \delta}f(x+v-\theta_1h)\int_{-\infty}^{+\infty}\chi_{-\frac{\delta}{h}+\theta_1\leq u\leq \frac{\delta}{h}+\theta_1}(u)\Big|K(u) - \frac{h^2}{2}\int s^2K(s)dsK''(u)\Big|\\
& \times \Bigg|\Bigg[K(\frac{x-\theta_1h-y+\theta_2h}{h}+u) - \frac{h^2}{2}\int s^2K(s)dsK''(\frac{x-\theta_1h-y+\theta_2h}{h}+u)\Bigg](h)^{-1}\Bigg|\, du
\end{split}
\end{equation*}
For every   $u$ 
\[\lim_{n \to +\infty}\Bigg|\Bigg[K(\frac{x-\theta_1h-y+\theta_2h}{h}+u) - \frac{h^2}{2}\int s^2K(s)dsK''(\frac{x-\theta_1h-y+\theta_2h}{h}+u)\Bigg](h)^{-1}\Bigg|=0.\]
Write 
\begin{equation*}
\begin{split}
&\Bigg|\Bigg[K(\frac{x-\theta_1h-y+\theta_2h}{h}+u) - \frac{h^2}{2}\int s^2K(s)dsK''(\frac{x-\theta_1h-y+\theta_2h}{h}+u)\Bigg](h)^{-1}\Bigg|=\\
&\Bigg|(\frac{x-\theta_1 h-y+\theta_2h}{h}+u)\Bigg[K(\frac{x-\theta_1h-y+\theta_2h}{h}+u) - \frac{h^2}{2}\int s^2K(s)dsK''(\frac{x-\theta_1h-y+\theta_2h}{h}+u)\Bigg]\Bigg|\\
& \times \Bigg|\frac{1}{x-\theta_1h-y+\theta_2h+hu}\Bigg|.
\end{split}
\end{equation*}
 We have 
\[\Big|\frac{1}{x-\theta_1h-y+\theta_2h+hu}\Big|=\frac{1}{|x-y||1-\frac{\theta_1-\theta_2-u}{x-y}h|}.\]
 Since  $|u|\leq \frac{\delta}{h}+\theta_1$ we may choose  $\delta$  small enough such that for  $n\geq n_0$ we have
 \[\Big|\frac{\theta_1-\theta_2-u}{x-y}h\Big|\leq \frac{3h+\delta}{|x-y|}=\eta<1.\]
 Therefore
 \begin{equation}
 \label{101}
 \Big|\frac{1}{x-\theta_1h-y+\theta_2h+hu}\Big|\leq \frac{1}{|x-y|(1-\eta)}
 \end{equation}
 since  $\mathbf{H_3}$ implies there exists $B$ such that 
 \[\Bigg|(\frac{x-\theta_1h-y+\theta_2h}{h}+u)\Bigg[K(\frac{x-\theta_1h-y+\theta_2h}{h}+u) - \frac{h^2}{2}\int s^2K(s)dsK''(\frac{x-\theta_1h-y+\theta_2h}{h}+u)\Bigg]\Bigg|\leq B\]
 Then, we have
\begin{equation*}
\begin{split}
\Bigg|K(\frac{x-\theta_1h-y+\theta_2h}{h}+u) - \frac{h^2}{2}\int s^2K(s)dsK''(\frac{x-\theta_1h-y+\theta_2h}{h}+u)\Bigg](h)^{-1}\Bigg|\leq \frac{B}{|x-y|(1-\eta)}
\end{split}
\end{equation*}
$\displaystyle |K(u)|$ being integrable, by dominated convergence
\[\int_{|v-\theta_1h|\leq \delta} \to 0 \quad \text{as} \quad n \to +\infty.\]
Let  $\displaystyle \int_{|v-\theta_1h|> \delta}$  write it in the form
\begin{equation*}
\begin{split}
&\int_{|v-\theta_1h|> \delta}=\int_{|v-\theta_1h|>\delta}\Bigg|v(h)^{-1}\Bigg[K(\frac{v}{h}) - \frac{h^2}{2}\int s^2K(s)dsK''(\frac{v}{h})\Bigg]\\
& \times \Bigg((h)^{-1}\Bigg[K(\frac{v+x-\theta_1-y+\theta_2h}{h}) - \frac{h^2}{2}\int s^2K(s)dsK''(\frac{v+x-\theta_1-y+\theta_2h}{h})\Bigg]\Bigg)\frac{f(x+v-\theta_1h)}{v}\Big|\, dv.
\end{split}
\end{equation*}
 We get
\begin{equation}
\begin{split}
\label{111}
\int_{|v-\theta_1h|> \delta}\leq &\frac{2}{\delta-\theta_1h}\sup_{|v-\theta_1h|> \delta}\Bigg|\frac{v}{h}\Bigg[K(\frac{v}{h}) - \frac{h^2}{2}\int s^2K(s)dsK''(\frac{v}{h})\Bigg]\Bigg|\int_{|v-\theta_1h|> \delta}((h)^{-1}\\
& \times \Bigg[K(\frac{v+x-\theta_1-y+\theta_2h}{h}) - \frac{h^2}{2}\int s^2K(s)dsK''(\frac{v+x-\theta_1-y+\theta_2h}{h})\Bigg])|f(x+v-\theta_1h)|\, dv. 
\end{split}
\end{equation}
Let the change of variable defined by
\[v+x-\theta_1h-y+\theta_2h= u.\]
Then
\begin{equation}
\begin{split}
\label{121}
\int_{|v-\theta_1h|> \delta}\leq \frac{2}{\delta-\theta_1h}&\sup_{|v-\theta_1h|> \delta}\Bigg|\frac{v}{h}\Bigg[K(\frac{v}{h}) - \frac{h^2}{2}\int s^2K(s)dsK''(\frac{v}{h}) \Bigg]\Bigg|\\
& \times \int_{\mathbb{R}}\Bigg|(h)^{-1}\Bigg[K(\frac{u}{h}) - \frac{h^2}{2}\int s^2K(s)dsK''(\frac{u}{h})\Bigg]\Bigg|f(u+y-\theta_2h)\, du. 
\end{split}
\end{equation}
Lemma \ref{lem211} (replacing  $\displaystyle K - \frac{h^2}{2}\int s^2K(s)dsK''$  by $\displaystyle \Bigg|K - \frac{h^2}{2}\int s^2K(s)dsK''\Bigg|$) and $(\mathbf{H_3})$ we have
\[|\int_{|v-\theta_1h|> \delta}| \to 0 \quad \text{as} \quad n \to +\infty\]
and the convergence is uniform.
Thus, the proof of lemma \ref{lem231}.
\end{preuve} 
\begin{rem}\rm
 \label{rem11}
 If condition $\mathbf{C_2}$ is verified, then the integral of the right hand-side of  \eqref{121} becomes
\begin{equation*}
\begin{split}
\int_{\mathbb{R}}\Bigg|(h)^{-1}\Bigg[K(\frac{u}{h}) - &\frac{h^2}{2}\int s^2K(s)dsK''(\frac{u}{h})\Bigg]\Bigg||f(u+y-\theta_2h)-  f(\frac{u}{h})+f(\frac{u}{h})|du\\
& \leq\int_{\mathbb{R}}\Bigg|(h)^{-1}\Bigg[K(\frac{u}{h}) - \frac{h^2}{2}\int s^2K(s)dsK''(\frac{u}{h})\Bigg]\Bigg|\int_{\frac{u}{h}}^{u+y-\theta_2h}|f^{'}(t)|\, dtdu\\
& \qquad +\int_{\mathbb{R}}f(\frac{u}{h})\Bigg|(h)^{-1}\Bigg[K(\frac{u}{h}) - \frac{h^2}{2}\int s^2K(s)dsK''(\frac{u}{h})\Bigg]\Bigg|\, du\\
&\leq \int_{\mathbb{R}}\Bigg|(h\lambda)^{-1}\Bigg[K(\frac{u}{h}) - \frac{h^2}{2}\int s^2K(s)dsK''(\frac{u}{h})\Bigg]\Bigg|\, du\int_{\mathbb{R}}|f^{'}(t)|\, dt+\int_{\mathbb{R}}|K(u) - \frac{h^2}{2}\int s^2K(s)dsK''(u)|f(u)\, du.
\end{split}
\end{equation*} 
The integrals of the right hand-side of this last inequality. Hence the theorem is valid under the hypothese  $\mathbf{C_2}$.
\end{rem}
\subsection*{Proof of theorem \ref{theo251}}
\begin{preuve}
We suppose condition $\mathbf{C_1}$ verified. Let $\displaystyle \Delta=[0, z]\times[0, z]$. We can write 
\begin{equation*}
\begin{split}
\sum_{0\leq i\neq j\leq[\frac{z}{h}]}\Big(1-\frac{ih}{z}\Big)^{\alpha}\Big(1-\frac{jh}{z}\Big)^{\alpha}\int_{\mathbb{R}}\Bigg|&\Bigg[K(\frac{u-ih}{h})  - \frac{h^2}{2}\int s^2K(s)dsK''(\frac{u-ih}{h})\Bigg]\\
&\times\Bigg[K(\frac{jh-u}{h}) - \frac{h^2}{2}\int s^2K(s)dsK''(\frac{jh-u}{h})\Bigg]\Bigg|f(u)\,du\\
& \qquad =\int_{\{(x,y) \in \Delta:|x-y|>0\}}\Phi_n(x,y)\, dxdy
\end{split}
\end{equation*} 
where 
\begin{equation*}
\begin{split}
\Phi_n(x,y)& =\frac{1}{(h)^2}\sum_{0\leq i\neq j\leq[\frac{z}{h}]}\chi_{\Delta_{h\,i}\times\Delta_{h,j}}(x,y)\Big(1-\frac{ih}{z}\Big)^{\alpha}\Big(1-\frac{jh}{z}\Big)^{\alpha}\\
& \quad \times \int_{\mathbb{R}}\Bigg|\Bigg[K(\frac{u-ih}{h}) - \frac{h^2}{2}\int s^2K(s)dsK''(\frac{u-ih}{h})\Bigg]\Bigg[K(\frac{jh-u}{h}) - \frac{h^2}{2}\int s^2K(s)dsK''(\frac{u-ih}{h})\Bigg]\Bigg|f(u)\, du.
\end{split}
\end{equation*}
If $\displaystyle (x,y)\in \Delta_{h,i}\times\Delta_{h,j} \quad i\neq j$  with the representation
\[x=hi+\theta_1h, \quad y=hj+\theta_2h \quad 0\leq \theta_l<1, \quad  l=1,2\]
We have
\begin{equation}
\begin{split}
\label{131}
\frac{1}{(h)^2}\Big(1-\frac{ih}{z}\Big)^{\alpha}&\Big(1-\frac{jh}{z}\Big)^{\alpha}\int_{\mathbb{R}}\Bigg|\Bigg[K(\frac{u-x+\theta_1h}{h}) - \frac{h^2}{2}\int s^2K(s)dsK''(\frac{u-x+\theta_1h}{h})\Bigg]\\
& \times \Bigg[K(\frac{u-y+\theta_2h}{h}) - \frac{h^2}{2}\int s^2K(s)dsK''(\frac{u-y+\theta_2h}{h})\Bigg]\Bigg|f(u)\, du\\
& \quad \leq \frac{1}{(h)^2}\int_{\mathbb{R}}\Bigg|\Bigg[K(\frac{u-x+\theta_1h}{h}) - \frac{h^2}{2}\int s^2K(s)dsK''(\frac{u-x+\theta_1h}{h})\Bigg]\\
& \quad \times \Bigg[K(\frac{u-y+\theta_2h}{h}) - \frac{h^2}{2}\int s^2K(s)dsK''(\frac{u-y+\theta_2h}{h})\Bigg]\Bigg|f(u)\, du
\end{split}
\end{equation}
The right hand-side of  \eqref{131} tends to  zero as  $\displaystyle n \to +\infty$  by lemma \ref{lem231}.\\
Let  $\displaystyle \delta=\frac{z}{2}$. Write
\begin{equation*}
\begin{split}
\frac{1}{(h)^2}\int_{\mathbb{R}}&\Bigg|\Bigg[K(\frac{u-x+\theta_1h}{h}) - \frac{h^2}{2}\int s^2K(s)dsK''(\frac{u-x+\theta_1h}{h})\Bigg]\\
&\times\Bigg[K(\frac{u-y+\theta_2h}{h}) - \frac{h^2}{2}\int s^2K(s)dsK''(\frac{u-y+\theta_2h}{h})\Bigg]\Bigg|f(u)\, du=\int_{|v|\leq\delta}+\int_{|v|>\delta}.
\end{split}
\end{equation*}
Then, we have
\begin{equation*}
\begin{split}
\int_{\{(x,y) \in \Delta:|x-y|>0\}}&\Phi_n(x,y)\, dxdy\\
&\leq\int_{\{(x,y) \in \Delta:|x-y|>0\}}\sum_{0\leq i\neq j\leq[\frac{z}{h}]}\chi_{\Delta_{h,i}\times\Delta_{h,j}}(x,y)(\int_{|v|\leq\delta}+\int_{|v|>\delta}).
\end{split}
\end{equation*}
The proof of the remainder is conducted as follow:\\
First consider  
\[\int_{\{(x,y) \in \Delta:|x-y|>0\}}\int_{|v|\leq\delta}.\]
Let $\displaystyle A=\sup_{x\in [0,z]}f(x)$. The notations being as in the proof lemma \ref{lem221}  with $\displaystyle \delta=\frac{z}{2}$, we have in accordance with inequality  \eqref{91}
\begin{equation*}
\begin{split} 
\int_{|v|\leq \delta}\leq A\int_{-\infty}^{+\infty}\chi_{-\frac{\delta}{h}\leq u\leq\frac{\delta}{h}}&\Bigg|K(u) - \frac{h^2}{2}\int s^2K(s)dsK''(u) \Bigg|\\
&\times\Bigg|\Bigg[K(\frac{x-\theta_1h-y+\theta_2h}{h}+u) - \frac{h^2}{2}\int s^2K(s)dsK''(\frac{x-\theta_1h-y+\theta_2h}{h}+u)\Bigg](h)^{-1}\Bigg|\, du.
\end{split}
\end{equation*} 
For every  $u$ 
\[\lim_{n \to +\infty}\Bigg|\Bigg[K(\frac{x-\theta_1h-y+\theta_2h}{h}+u) - \frac{h^2}{2}\int s^2K(s)dsK''(\frac{x-\theta_1h-y+\theta_2h}{h}+u)\Bigg](h)^{-1}\Big|=0.\]
We have 
\begin{equation*}
\begin{split}
&\Bigg|\Bigg[K(\frac{x-\theta_1h-y+\theta_2h}{h}+u) - \frac{h^2}{2}\int s^2K(s)dsK''(\frac{x-\theta_1h-y+\theta_2h}{h}+u)\Bigg](h)^{-1}\Bigg|\\
&=\Bigg|\Bigg[K(\frac{x-\theta_1h-y+\theta_2h}{h}+u) - \frac{h^2}{2}\int s^2K(s)dsK''(\frac{x-\theta_1h-y+\theta_2h}{h}+u)\Bigg]\\
&- \Bigg[K(\frac{2z+x-\theta_1-y+\theta_2h}{h}+u) - \frac{h^2}{2}\int s^2K(s)dsK''(\frac{2z+x-\theta_1-y+\theta_2h}{h}+u)\Bigg]\\
&+\Bigg[K(\frac{2z+x-\theta_1-y+\theta_2h}{h}+u) - \frac{h^2}{2}\int s^2K(s)dsK''(\frac{2z+x-\theta_1-y+\theta_2h}{h}+u)\Bigg]\Bigg|(h)^{-1}\\
&\leq \Bigg(\lambda(\frac{2z}{h})+\Bigg|K(\frac{2z+x-\theta_1h-y+\theta_2h}{h}+u) - \frac{h^2}{2}\int s^2K(s)dsK''(\frac{2z+x-\theta_1-y+\theta_2h}{h}+u)\Bigg|\Bigg)(h)^{-1}.
\end{split}
\end{equation*}
Moreover,
\begin{equation*}
\begin{split}
&\Bigg|K(\frac{2z+x-\theta_1-y+\theta_2h}{h}+u) - \frac{h^2}{2}\int s^2K(s)dsK''(\frac{2z+x-\theta_1-y+\theta_2h}{h}+u)\Bigg|(h)^{-1}\\
&=\Big|\frac{2z+x-y+hu}{h}\Big|\Bigg|K(\frac{2z+x-\theta_1-y+\theta_2h}{h}+u) - \frac{h^2}{2}\int s^2K(s)dsK''(\frac{2z+x-\theta_1-y+\theta_2h}{h}+u)\Bigg|\frac{1}{|2z+x-y+hu|}  
\end{split}
\end{equation*}
Let $\displaystyle B=\sup_{y\in\mathbb{R}}|y||K(y) - \frac{h^2}{2}\int s^2K(s)dsK''(y)|$ and  $\displaystyle C =\sup_{y\in\mathbb{R}}|K(y) - \frac{h^2}{2}\int s^2K(s)dsK''(y)|$, then we have
\begin{equation*}
\begin{split}
\Big|\frac{2z+x-y+hu}{h}\Big|\Bigg|K(\frac{2z+x-\theta_1-y+\theta_2h}{h}+u) - \frac{h^2}{2}\int s^2K(s)dsK''(\frac{2z+x-\theta_1-y+\theta_2h}{h}+u)\Bigg|\leq B+2hC. 
\end{split}
\end{equation*}
Therefore, we have 
\begin{equation*}
\begin{split}
\Bigg|K(\frac{2z+x-\theta_1-y+\theta_2h}{h}+u) - \frac{h^2}{2}\int s^2K(s)dsK''(\frac{2z+x-\theta_1-y+\theta_2h}{h}+u)\Bigg|(h)^{-1} \leq \frac{B+2hC}{|2z+x-y+hu|}
\end{split}
\end{equation*}
Hence
\begin{equation*}
\begin{split}
\Bigg|\Bigg[K(\frac{x-\theta_1h-y+\theta_2h}{h}+u) - \frac{h^2}{2}\int s^2K(s)dsK''(\frac{x-\theta_1h-y+\theta_2h}{h}+u)\Bigg](h)^{-1}\Bigg|\leq \lambda(\frac{2z}{h})+\frac{B+2hC}{|2z+x-y+hu|}
\end{split}
\end{equation*}
We conclude that for $h$ small enough 
\begin{equation*}
\begin{split}
\int_{|v|\leq \delta} &\leq \frac{A}{|2z+x-y+hu|}\int_{\mathbb{R}}\Bigg|K(u) - \frac{h^2}{2}\int s^2K(s)dsK''(u)\Bigg|(B+2hC)\, du< \frac{AD}{|2z+x-y+hu|}\\
&\leq \frac{AD}{|2z+x-y+hu|} \\
&\leq \frac{AD}{(2z+x-y+hu)}
\end{split}
\end{equation*}
$\displaystyle D$ being the finite bound of $\displaystyle \int_{\mathbb{R}}\Bigg|K(u) - \frac{h^2}{2}\int s^2K(s)dsK''(u)\Bigg|(B+2C) du$.\\
Finally, we have
\[\int_{|v|\leq \delta} \leq \frac{AD}{(2z+x-y+hu)}+O(h).\]
Since  $\displaystyle -\delta\leq hu\leq\delta$ \, we have  $\displaystyle \frac{z}{2}\leq 2z+x-y+hu\leq \frac{7z}{2}$. \\
Hence
\[\int_{|v|\leq \delta} \leq \frac{2AD}{z}+O(h).\]
Therefore, by Lebesgue-dominated convergence, we have
\begin{equation}
\begin{split}
\label{141}
\lim_{n \to +\infty}\int\int_{\Delta}\Big(\int_{|v|\leq \delta}\Big)\, dxdy=\int\int_{\Delta}\lim_{n \to +\infty}\Big(\int_{|v|\leq \delta}\Big)\, dxdy=0
\end{split}
\end{equation}
Consider then $\displaystyle \int_{|v| > \delta}$.\\
We use the second part, by analogous reasoning, of the proof of lemma \ref{lem221} 
\begin{equation}
\begin{split}
\label{151}
\int_{|v| > \delta}\leq \frac{2}{\delta}\sup_{|v|> \delta}\Big|\frac{v}{h}\Big|\Bigg|K(\frac{v}{h}) - \frac{h^2}{2}\int s^2K(s)dsK''(\frac{v}{h}) \Bigg|\int_{\mathbb{R}}\Big|(h)^{-1}\Bigg[K(\frac{u}{h}) - \frac{h^2}{2}\int s^2K(s)dsK''(\frac{u}{h})\Bigg]\Big|f(u+y-\theta_2h)\, du. 
\end{split}
\end{equation}
We have
\[\int_{|v| > \delta} \to 0, \qquad n \to +\infty \quad \text{uniformly}.\]
Hence
\[\lim_{n \to +\infty} \int_{\Delta}\int_{\mathbb{R}} \to 0 , \quad n \to +\infty\]
since $\displaystyle \Delta$ is bounded. The proof of the lemma is complete.
\end{preuve}
\begin{rem}
\label{rem21}
If  $\mathbf{C_2}$ is verified, then the theorem is a gain valid. 
Indeed it suffices to apply  remark \ref{rem11} to inequality  \eqref{151}.
\end{rem}
\subsection*{Proof of  theorem \ref{theo231}}
\begin{preuve}
We suppose condition $(\mathbf{C_1})$ satisfied.
\begin{equation}
\begin{split}
\label{161}
n\mathbb{V}ar(P_n(z,\alpha)) &=\mathbb{E}\Bigg(\sum_{i=0}^{[\frac{z}{h}]}\Big(1-\frac{ih}{z}\Big)^{\alpha}\Bigg[K\Big(\frac{ih-X_j}{h}\Big) - \frac{h^2}{2}\int s^2K(s)dsK''\Big(\frac{ih-X_j}{h}\Big)\Bigg]\Bigg)^2\\
& \quad -\mathbb{E}^2\Bigg(\sum_{i=0}^{[\frac{z}{h}]}\Big(1-\frac{ih}{z}\Big)^{\alpha}\Bigg[K\Big(\frac{ih-X_j}{h}\Big) - \frac{h^2}{2}\int s^2K(s)dsK''\Big(\frac{ih-X_j}{h}\Big)\Bigg]\Bigg)
\end{split}
\end{equation}
Since
\begin{equation}
\begin{split}
\label{171}
&\mathbb{E}\Bigg(\sum_{i=0}^{[\frac{z}{h}]}\Big(1-\frac{ih}{z}\Big)^{\alpha}\Bigg[K\Big(\frac{ih-X_j}{h}\Big) - \frac{h^2}{2}\int s^2K(s)dsK''\Big(\frac{ih-X_j}{h}\Big) \Bigg]\Bigg)^2\\
&=\mathbb{E}\Bigg[\Bigg(\sum_{i=0}^{[\frac{z}{h}]}\Big(1-\frac{ih}{z}\Big)^{\alpha}\Bigg[K\Big(\frac{ih-X_j}{h}\Big) - \frac{h^2}{2}\int s^2K(s)dsK''\Big(\frac{ih-X_j}{h}\Big)\Bigg]\Bigg)^2\Bigg]\\
&=\mathbb{E}\Bigg[\bigg\{\sum_{i=0}^{[\frac{z}{h}]}\Big(1-\frac{ih}{z}\Big)^{2\alpha}\Bigg[K\Big(\frac{X_k-ih}{h}\Big) - \frac{h^2}{2}\int s^2K(s)dsK''\Big(\frac{X_k-ih}{h}\Big)\Bigg]^2\\
& \quad +\sum_{i\neq j}^{[\frac{z}{h}]}\Big(1-\frac{ih}{z}\Big)^{\alpha}\Big(1-\frac{jh}{z}\Big)^{\alpha}\Bigg[K\Big(\frac{X_k-ih}{h}\Big) - \frac{h^2}{2}\int s^2K(s)dsK''\Big(\frac{X_k-ih}{h}\Big)\Bigg]\\
& \qquad \times \Bigg[K\Big(\frac{X_k-jh}{h}\Big) - \frac{h^2}{2}\int s^2K(s)dsK''\Big(\frac{X_k-jh}{h}\Big)\Bigg]\bigg\}\Bigg]\\
&=\sum_{i=0}^{[\frac{z}{h}]}\Big(1-\frac{ih}{z}\Big)^{2\alpha}\int_{\mathbb{R}}\Bigg[K\Big(\frac{u-ih}{h}\Big) - \frac{h^2}{2}\int s^2K(s)dsK''\Big(\frac{u-ih}{h}\Big)\Bigg]^2f(u)\, du\\
&+\sum_{i\neq j}^{[\frac{z}{h}]}\Big(1-\frac{ih}{z}\Big)^{\alpha}\Big(1-\frac{jh}{z}\Big)^{\alpha}\int_{\mathbb{R}}\Bigg[K\Big(\frac{u-ih}{h}\Big) - \frac{h^2}{2}\int s^2K(s)dsK''\Big(\frac{u-ih}{h}\Big)\Bigg]\\
& \qquad \times \Bigg[K\Big(\frac{jh-u}{h}\Big) - \frac{h^2}{2}\int s^2K(s)dsK''\Big(\frac{jh-u}{h}\Big)\Bigg]f(u)\, du.
\end{split}
\end{equation}
It follows by corollary \ref{cor211}, of  lemma \ref{lem221} and theorem \ref{theo251} that, as  $\displaystyle n \to +\infty$, the first term of the right hand-side of \eqref{171} tends to $\displaystyle \Bigg(\int_{\mathbb{R}}K^2(y)P(z,2\alpha)\Bigg)$ and the second term tends to zero uniformly on $\displaystyle [0, b]$.\\
then 
\[n\mathbb{V}ar(P_n(z,\alpha)) \to \Bigg(\int_{\mathbb{R}}K^2(y) P(z,2\alpha)-P^2(z,\alpha)\Bigg)\]
According to the Lemma \ref{lem211}, we have 
\[\mathbb{E}(P_n(z,\alpha)) \to P(z,\alpha) \quad n \to +\infty\]
therefore
\[\mathbb{E}^2(P_n(z,\alpha)) \to P^2(z,\alpha) \quad n \to +\infty.\]
Define
\[biais(P_n(z,\alpha))=\mathbb{E}(P_n(z,\alpha))-\mathbb{V}ar(P_n(z,\alpha)).\]
We have 
\[\mathbb{E}(P_n(z,\alpha)-P(z,\alpha))^2=biais^2(P_n(z,\alpha))+\mathbb{V}ar(P_n(z,\alpha))\]
and
\[\Bigg|\Bigg(\int_{\mathbb{R}}K^2(y)P(z,2\alpha)-P^2(z,\alpha)\Bigg)\Bigg|\leq \int_{\mathbb{R}}K^2(y) + 1.\]
Hence
\[\mathbb{V}ar(P_n(z,\alpha)) = O(\frac{1}{n}).\]
By Lemma \ref{lem211}, we have 
\[|\mathbb{E}(P_n(z,\alpha))-P_n(z,\alpha)| \to 0 , \quad n \to +\infty,\]
therefore
\[biais^2(P_n(z,\alpha)) \to 0, \quad n \to +\infty\]
hence 
\[\mathbb{E}(P_n(z,\alpha)-P(z,\alpha))^2  \to 0, \quad n \to +\infty.\]
If  condition $\mathbf{C_2}$ is satisfied, the Theorem is again valid, by  Corollary \ref{cor211}, of Lemma \ref{lem221}  and using Remark \ref{rem21} of Theorem \ref{theo251}  and  Theorem \ref{theo221}. 
\end{preuve}

\end{document}